\documentclass[11pt,a4paper]{article}

\usepackage[T1]{fontenc}
\usepackage{lmodern}
\usepackage[hmargin=23mm,
vmargin=25mm]{geometry}
\usepackage{amsmath,amssymb,amsthm}
\usepackage{microtype}
\usepackage{tikz}
\usepackage[font=small,labelfont=bf,labelsep=period,skip=6pt]{caption}
\usepackage[hidelinks]{hyperref}
\usetikzlibrary{arrows.meta}

\usepackage{authblk}

\usepackage{stmaryrd}

\usepackage{placeins}

\newcommand \brk[1]
{
	\llbracket #1\rrbracket
}

\hypersetup{pdftitle={Spectral Radii of Connected Graphs}}
\newtheorem{theorem}
{Theorem}[section]
\newtheorem{lemma}[theorem]{Lemma}
\newtheorem{corollary}[theorem]{Corollary}
\theoremstyle{definition}
\newtheorem{example}
{Example}[section]

\newtheorem{prob}{Problem}

\newcommand{\diag}{\operatorname{diag}}
\newcommand{\trans}{\mathsf T}
\tikzset{
  graph/.style={line width=.68pt,line cap=round,line join=round,font=\small},
  vertex/.style={circle,fill=black,inner sep=0pt,minimum size=4.4pt},
  branchvertex/.style={vertex,fill=white,draw=black,minimum size=4.7pt},
  branch/.style={line width=.68pt,draw=black!65}
}

\usepackage{enumerate}

\usepackage{mathrsfs}

\def \scrs {{\mathscr S}}
\def \scrq {{\mathscr Q}}
\newcommand{\setw}{{\mathscr W}}

\newcommand{\scrt}{{\mathscr T}}

\def \wisetwra 
{\widehat \setw_{r}(\alpha)}

\def \bfd {\mathbf d}
\def \bfb {\mathbf b}
\def \bfa {\mathbf a}

\def\sgap{\vspace{3mm}}

\newcommand{\equ}[2]
{
\begin{equation}
\label{#1}#2
\end{equation}
}

\begin{document}

	\title{Spectral radii of connected graphs with prescribed maximum degree
	}

\author[]{Fengming Dong\thanks{
Email: fengming.dong@nie.edu.sg and donggraph@163.com}}
		
		\affil[]{National Institute of Education, Nanyang Technological University, Singapore}
		
		\date{}
		
		\maketitle

\begin{abstract}
	For a real number $\alpha>0$ and a positive integer $r$, we give recursive
	characterizations of the existence of finite simple connected graphs $G$
	with adjacency spectral radius $\rho(G)=\alpha$ and maximum degree 
	$\Delta(G)\le r$ (or 
	$=r$).
Adding one vertex 
to a graph $F$ at a time
to update 
	the Schur-complement 
	$q=\alpha-\bfb^{\mathsf T}B^{-1}\bfb$, where
	$B=\alpha I-A(F)$ 
	is positive definite and
	$\bfb$ is 
	 a nonzero
	$0$--$1$ that specifies the neighbours of $z$ in $F$. 
	Realizability is
	characterized by admissible constructions with positive intermediate
	outputs $q$ and a final zero output. We also establish a reduction for a
	graph $G$ obtained by attaching pairwise disjoint rooted trees, each by
	a single edge, to a nonempty connected graph $H$. Assume that
	$\alpha I-A(T)\succ0$ for every attached tree $T$. For each
	$z\in V(H)$, let $q_{z,1},\ldots,q_{z,s_z}$ be the positive root pivots
	obtained by the scalar tree recursion for the trees attached to $z$, and
	set
	\[
	a_z=\alpha-\sum_{j=1}^{s_z}q_{z,j}^{-1},
	\qquad
	M_H(\mathbf a)=\operatorname{diag}(a_z:z\in V(H))-A(H),
	\]
	where an empty sum is zero. Then
	$\rho(G)=\alpha$
	if and only if 
$M_H(\mathbf a)\succeq0$
and $\det M_H(\mathbf a)=0$.
	Taking $H=K_1$ recovers the tree
	criteria, while $\delta(H)\geq2$ gives criteria for graphs with
	prescribed 2-core $H$. In particular, for unicyclic graphs of prescribed
	cycle length, path continuants reduce the spectral condition to
	positivity inequalities and one scalar equation.
\end{abstract}




\section{Introduction}
\label{se:intro}

All graphs considered are finite, simple, and undirected. For a graph $G$,
let $A(G)$ denote its adjacency matrix, $\rho(G)$ its adjacency spectral
radius, and $\Delta(G)$ its maximum degree. Since $A(G)$ is symmetric and
nonnegative, $\rho(G)$ is its largest eigenvalue. A graph is called
nontrivial if it has at least two vertices. For each positive integer $r$,
let
\begin{equation}\label{eq:intro-Q}
	\mathcal Q_r=
	\{\rho(G):G\text{ is nontrivial and connected},\ \Delta(G)\leq r\},
	\qquad
	\mathcal Q=\bigcup_{r\geq1}\mathcal Q_r.
\end{equation}
For a real target $\alpha>0$, we study the existence of graphs with
spectral radius $\alpha$ and maximum degree at most $r$ or exactly $r$.
The standard inequalities
\[
\sqrt{\Delta(G)}\leq\rho(G)\leq\Delta(G)
\]
give the necessary numerical restriction
\[
\lceil\alpha\rceil\leq r\leq\lfloor\alpha^2\rfloor
\qquad\text{when }\rho(G)=\alpha\text{ and }\Delta(G)=r;
\]
see \cite[Proposition 3.1.2 and Section 3.15, Exercise 1]{BrouwerHaemers}.
The first two set $\scrq_1$ and $\scrq_2$ are explicitly given below:
\[
\mathcal Q_1=\{1\},\qquad
\mathcal Q_2=
\left\{2\cos\frac{\pi}{n+1}:n\geq2\right\}\cup\{2\}.
\]
Indeed, a nontrivial connected graph of maximum degree at most two is
a path or a cycle, with
$\rho(P_n)=2\cos(\pi/(n+1))$ and $\rho(C_m)=2$; see
\cite[Sections 1.4.3--1.4.4]{BrouwerHaemers}.

Determining the elements of $\mathcal Q$ is an inverse spectral problem.
If $\alpha\in\mathcal Q$, then $\alpha$ is a totally real algebraic
integer and every algebraic conjugate $\beta$ of $\alpha$ satisfies
$|\beta|\leq\alpha$. These conditions are not sufficient in general;
see \cite[Remark 6.2]{Schildkraut}. Classical results of Smith
\cite{Smith} describe the connected graphs with spectral radius at most
two. Ding, Ge, and Liu \cite{DingGeLiu} characterized realizability
among quadratic algebraic integers and among positive real numbers at
most two. Prescribing a maximum degree or a nonempty 2-core adds
structural requirements to this realizability problem.

For trees, Dong and Zhang \cite[Theorems 2 and 3]{DongZhang} introduced
the scalar recursion
\[
q=\alpha-\sum_{i=1}^{s}q_i^{-1}.
\]
Their set $\mathscr W_r(\alpha)$ is the smallest subset of
$[0,\infty)$ containing $\alpha$ and admitting this output 
$q$
from positive
inputs $q_i\in\mathscr W_r(\alpha)$ whenever either
$1\leq s<r$ and $q\geq0$, or $s=r$ and $q=0$. Repetition of inputs
is allowed. Put
$\mathscr W_r^+(\alpha)=\mathscr W_r(\alpha)\setminus\{0\}$.
Their results give
\[
0\in\mathscr W_r(\alpha)
\quad\Longleftrightarrow\quad
\text{there exists a tree }T\text{ with }
\rho(T)=\alpha\text{ and }\Delta(T)\leq r.
\]
The exact-degree condition is
\[
\begin{gathered}
	\exists\, \text{ a tree }
	T\text{ with }
	\rho(T)=\alpha\text{ and }\Delta(T)=r\quad
	\Longleftrightarrow\quad
	\exists\, 
	q_1,\ldots,q_r\in\mathscr W_r^+(\alpha)
	\text{ such that }\sum_{i=1}^{r}q_i^{-1}=\alpha.
\end{gathered}
\]
Positive outputs represent rooted branches whose vertices have at most
$r-1$ children; a final zero output allows up to $r$ branches at the
last root.

Our first characterization extends this construction to connected
graphs by retaining a matrix and its degree tuple. Start with
$B=[\alpha]$ and $\mathbf d=(0)$, corresponding to $K_1$. At a step
with $B\succ0$ of order $n$, choose $b\in\{0,1\}^n$ satisfying
\[
1\leq s:=\sum_{i=1}^{n}b_i\leq r,
\qquad d_i+b_i\leq r\quad(1\leq i\leq n),
\]
and form
\[
B^+=\begin{pmatrix}B&-b\\-b^{\mathsf T}&\alpha\end{pmatrix},
\qquad
\mathbf d^+=(d_1+b_1,\ldots,d_n+b_n,s).
\]
The next Schur-complement output is
\[
q=\alpha-b^{\mathsf T}B^{-1}b
=\frac{\det B^+}{\det B}.
\]
Retain the update if $q\geq0$, continue from it if $q>0$, and
terminate it if $q=0$. Let $\widehat{\mathscr W}_r(\alpha)$ consist
of $\alpha$ and all nonnegative outputs obtained by finite sequences
of these updates, as formalized in Section~\ref{se2}. Theorem~\ref{th:main} proves
\[
\alpha\in\mathcal Q_r
\quad\Longleftrightarrow\quad
0\in\widehat{\mathscr W}_r(\alpha).
\]
Maximum degree exactly $r$ is obtained precisely when a terminal update
has $\max_{1\leq i\leq n+1}d_i^+=r$. Since
$\Delta(G)\leq\rho(G)^2$, Corollary~\ref{cor:all} further yields
\[
\mathcal Q=
\left\{\alpha\geq1:
0\in\widehat{\mathscr W}_{\lfloor\alpha^2\rfloor}(\alpha)\right\}.
\]

Our second characterization reduces the matrix size when the part
supporting the cycles is prescribed. Let $H$ be a nonempty connected
graph, and form $G$ by attaching pairwise disjoint rooted trees to its
vertices, each by a single edge. Assume that
$M_\alpha(T):=\alpha I-A(T)\succ0$ for every attached tree $T$.
Eliminating children before parents gives positive pivots satisfying
the scalar recursion. If the root labels of the $s_z$ trees attached
to $z\in V(H)$ are $q_{z,1},\ldots,q_{z,s_z}$, put
\[
a_z=\alpha-\sum_{j=1}^{s_z}q_{z,j}^{-1},
\qquad
M_H(\mathbf a)=\operatorname{diag}(a_z:z\in V(H))-A(H),
\]
where an empty sum is zero. The elimination gives the congruence
\[
\alpha I-A(G)\ \text{congruent to}\
\operatorname{diag}(D,M_H(\mathbf a)),\qquad D\succ0,
\]
with the block $D$ omitted when there are no attached trees.
Consequently, 
Theorem~\ref{se3-th1} and Corollary~\ref{se3-th1-c0} give
\[
\begin{aligned}
	\rho(G)=\alpha
	&\quad\Longleftrightarrow\quad
	M_H(\mathbf a)\succeq0\ \text{and}\ \det M_H(\mathbf a)=0,\\
	\rho(G)<\alpha
	&\quad\Longleftrightarrow\quad
	M_H(\mathbf a)\succ0.
\end{aligned}
\]

For a prescribed degree bound $r\geq\max\{1,\Delta(H)\}$, choose
the attached rooted trees from the positive constructions of the
scalar recursion in 
Section~\ref{sec3-1}. Every attached vertex then has at
most $r-1$ children, and every root label belongs to
$\mathscr W_r^+(\alpha)$. The remaining degree bounds are exactly
\[
0\leq s_z\leq r-d_H(z)\qquad(z\in V(H)).
\]
Under these bounds, maximum degree exactly $r$ is attained if and
only if
\[
\begin{gathered}
	d_H(z)+s_z=r\quad\text{for some }z\in V(H),\\
	\text{or}\qquad
	c(v)=r-1\quad\text{for some }v\in V(G)\setminus V(H),
\end{gathered}
\]
where $c(v)$ is the number of children of $v$ in its attached rooted
tree. In the existence criteria, a branch label must be realized by
a rooted tree witnessing the second alternative when that alternative
is used. Thus the exact-degree condition records where the degree is
attained, as well as the numerical values of the root labels.

The 2-core $G'$ of a connected graph $G$ is obtained by repeatedly
deleting vertices of degree at most one. It is empty for a tree and
is nonempty and connected when $G$ contains a cycle. Every component
outside a nonempty connected induced subgraph $H$ containing $G'$
is a tree joined to $H$ by a single edge. Corollary~\ref{th:core} therefore
turns the preceding reduction into an existence criterion for such
extensions of $H$. If $\delta(H)\geq2$, the constructed graph has
$G'=H$, and the bounded-degree criterion becomes
\[
\begin{gathered}
	\text{there exists a connected }G\text{ with }
	G'\cong H,\ \rho(G)=\alpha,\ \Delta(G)\leq r\\
	\Longleftrightarrow\quad
	\text{there are admissible branch choices with }
	M_H(\mathbf a)\succeq0\text{ and }\det M_H(\mathbf a)=0.
\end{gathered}
\]
Here admissibility means the bounds on $s_z$ above and positive root
labels from $\mathscr W_r^+(\alpha)$. The specialization $H=K_1$
reduces the matrix condition to
$\alpha-\sum_{j=1}^{s_z}q_{z,j}^{-1}=0$, recovering the tree criteria.

For a unicyclic graph, the 2-core is its unique cycle. Fix $m\geq3$
and $r\geq2$, label the cycle vertices $v_1,\ldots,v_m$ in cyclic
order, and write $s_i=s_{v_i}$ and $a_i=a_{v_i}$, so
$0\leq s_i\leq r-2$. For $1\leq i\leq j\leq m-1$, let $D(i,j)$
be the determinant of the tridiagonal path matrix with diagonal
$a_i,\ldots,a_j$ and adjacent off-diagonal entries $-1$. These
continuants satisfy
\[
\begin{aligned}
D(i,i-1)&=1,\qquad D(i,i)=a_i,
\qquad
D(i,j)=a_jD(i,j-1)-D(i,j-2)\qquad(i<j).
\end{aligned}
\]
Theorem 4.2 shows that a unicyclic graph of cycle length $m$,
spectral radius $\alpha$, and maximum degree at most $r$ exists if
and only if admissible branch choices can be made such that
\[
\begin{gathered}
	D(1,t)>0\quad(1\leq t\leq m-1),\qquad
	a_mD(1,m-1)=D(1,m-2)+D(2,m-1)+2.
\end{gathered}
\]
For maximum degree exactly $r$, the chosen realization must additionally
have $s_i=r-2$ at a cycle vertex or a vertex with $r-1$ children in an
attached tree. In particular, when $r=2$, all $s_i$ vanish and the only
realization of the prescribed cycle length is $C_m$, with spectral
radius $2$.

Schur-complement congruences and their relation to inertia are standard;
see \cite{HornJohnson,ZhangSchur}. Related elimination procedures
underlie the eigenvalue-location algorithms of Jacobs and Trevisan
\cite{JacobsTrevisan} for trees and Braga, Rodrigues, and Trevisan
\cite{BragaRodriguesTrevisan} for unicyclic graphs. Here the calculations
give recursive existence criteria with explicit degree restrictions and
a prescribed base graph or 2-core. For algebraic $\alpha$, a finite
construction provides an exactly checkable certificate. No stopping
bound on negative instances is established. Section 5 asks about
effective recognition, square-root spectral radii of unicyclic graphs,
infinitely many realizations with one fixed 2-core, and accumulation
points for a fixed cycle. The latter questions are related to the
equal-radius constructions of 
Seeger and Sossa \cite{Seeger} and the Hoffman program surveyed
in \cite{Wang}.

\section{A matrix recursion for connected graphs}

\label{se2}

\subsection{Matrix preliminaries}

For a real number \(\alpha>0\) and a simple graph $G$, set \(M_\alpha(G)=\alpha I_n-A(G)\), 
where $n=|V(G)|$.
Write \(M\succ0\) and
\(M\succeq0\) for positive definiteness and positive
semidefiniteness, respectively. For background on adjacency spectra,
see~\cite[Chapter~1]{BrouwerHaemers}. 
The spectral theorem gives
\begin{align}
 M_\alpha(G)\succ0
   &\ \Longleftrightarrow\ \alpha>\rho(G),\label{eq:pd}\\
 M_\alpha(G)\succeq0\text{ and singular}
   &\ \Longleftrightarrow\ \alpha=\rho(G).\label{eq:psd}
\end{align}
If \(G\) is connected and \(\alpha=\rho(G)\), the Perron--Frobenius
theorem gives a strictly positive vector spanning \(\ker M_\alpha(G)\).
The following lemma is an immediate consequence of eigenvalue interlacing and its equality case; see 
\cite[Theorem 2.5.1(i),(ii)]{BrouwerHaemers}.


\begin{lemma}
	[\cite{BrouwerHaemers}]
	\label{le:prin}
Let \(M\) be a real symmetric positive semidefinite matrix whose kernel
is spanned by a vector with every coordinate nonzero. Every nonempty
proper principal submatrix of \(M\) is positive definite.
\end{lemma}

We will repeatedly use the Schur complement of a positive definite
block. 

Two real symmetric matrices \(M\) and \(N\) of the same order are congruent if there exists an invertible real matrix \(P\) such that
$
N=P^{\mathsf T}MP.
$

\begin{lemma}\label{se2-le1}
	Let \(B\), 
	\(F\)  and $U$ be 
	$n\times n$, 
	 $m\times m$
	and $n\times m$
	real matrices, respectively.
If both $B$ and \(F\) are symmetric and $B\succ 0$,  then
\begin{equation}
	\label{eq:congr}
M:= 
\begin{pmatrix}
B&-U\\-U^{\trans}&F
\end{pmatrix}
 \ \text{is congruent to}\
 \diag\bigl(B,F-U^{\trans}
 B^{-1}U\bigr)
\end{equation}
and 
\equ{le:prin-e2}
{
\det M=\det B \det (F-U^{\trans}
B^{-1}U).
}
\end{lemma}

\begin{proof}
Write \(S=F-U^{\mathsf T}B^{-1}U\) for the Schur complement of \(B\).
To prove \eqref{eq:congr}, take
\[
P=\begin{pmatrix}I_n&B^{-1}U\\0&I_m\end{pmatrix}.
\]
The matrix \(P\) is invertible with \(\det P=1\). Since \(B^{-1}\) is symmetric, the two block multiplications give
\begin{align}
	MP
	&=\begin{pmatrix}
		B&BB^{-1}U-U\\
		-U^{\mathsf T}&F-U^{\mathsf T}B^{-1}U
	\end{pmatrix}
	=\begin{pmatrix}B&0\\-U^{\mathsf T}&S\end{pmatrix},
	\\[4pt]
	P^{\mathsf T}MP
	&=\begin{pmatrix}
	I_n&0\\U^{\mathsf T}B^{-1}&I_m\end{pmatrix}
	\begin{pmatrix}B&0\\-U^{\mathsf T}&S\end{pmatrix}
	=\begin{pmatrix}B&0\\0&S\end{pmatrix}.
	\label{eq:cong2}
\end{align}
This is the congruence asserted in  \eqref{eq:congr}, 
and \eqref{le:prin-e2}
follows from \eqref{eq:cong2}.
\end{proof}

 Lemma~\ref{se2-le1}
 gives the following corollary.

\begin{corollary}
	\label{se2-co1}
	Let \(B\) be a real symmetric  \(n\times n\) matrix,
	 let \(U\in\mathbb R^{n\times1}\), and let \(a>0\). If $B\succ 0$, then the following statements are equivalent.
\begin{enumerate}[(i)]
	\item $\begin{pmatrix}
		B&-U\\-U^{\trans}&a
	\end{pmatrix}$ is 
positive definite (resp., positive semidefinite; positive semidefinite and singular);

\item $a-U^{\trans} B^{-1} U>0$
(resp., \(\ge0\); \(=0\)).

\item $B-\frac 1a 
UU^{\trans}$
 is 
positive definite (resp., positive semidefinite; positive semidefinite and singular).
\end{enumerate}
\end{corollary}

\begin{proof}
	Write
	\[
	M=\begin{pmatrix}B&-U\\-U^{\mathsf T}&a\end{pmatrix},
	\qquad
	q=a-U^{\mathsf T}B^{-1}U,
	\qquad
	C=B-\frac1aUU^{\mathsf T}.
	\]
	Applying Lemma~\ref{se2-le1} to \(B\), and then to \(a\) after interchanging the blocks, gives
	\[
	M\ \text{congruent to }\operatorname{diag}(B,q),
	\qquad
	M\ \text{congruent to }\operatorname{diag}(a,C).
	\]
	Since \(B\succ0\) and \(a>0\), the three equivalences follow. In particular, \(q=0\) holds exactly when either matrix is singular, and both matrices are then positive semidefinite.
\end{proof}

\subsection{Construction and characterization}

Fix \(\alpha>0\) and an integer \(r\ge1\). A state \((B,\mathbf d,H)\) consists of a connected simple graph \(H\), with ordered vertices \(v_1,\ldots,v_n\) and \(\Delta(H)\le r\), the matrix \(B=M_\alpha(H)\succeq0\), and its degree tuple \(\mathbf d=(d_H(v_1),\ldots,d_H(v_n))\). A state \((B,\mathbf d,H)\)  is called nonterminal if \(B\succ0\), and terminal if \(B\) is singular.
All states with respect to $\alpha$ and $r$ 
are generated by 
the following rules.
\begin{enumerate}[(i)]
	\item The initial state is
$(B,\mathbf d, H)=([\alpha],(0), K_1),
$
where \([\alpha]\) is the \(1\times1\) matrix with entry \(\alpha\),
namely \(M_\alpha(K_1)\).

\item 
Let \((B,\bfd, H)\) be a 
nonterminal state   corresponding to
a connected simple graph $H$
with ordered vertices  $v_1,\ldots,v_n$,
where $\bfd=(d_1,\ldots,d_n)$
with $d_i=d_H(v_i)$.
Choose any column
\(\bfb=(b_1,\ldots,b_n)^{\trans}\in\{0,1\}^n\) satisfying
\begin{equation}
	\label{eq:admis}
 1\le s:=\sum_{i=1}^n b_i\le r,\qquad d_i+b_i\le r
 \quad(i\in \brk{n}),
\end{equation}
where $\brk{n}:=\{1,\ldots,n\}$
and $\brk{n}=\emptyset$ if $n<1$.
Let \(H^+\) be the graph 
obtained from $H$ by 
adding a new vertex \(v_{n+1}\)
and new edges joining 
$v_{n+1}$ to vertices in $H$
such that the neighbours of $v_{n+1}$ in $H^+$
are exactly the
vertices \(v_i\) for which \(b_i=1\).
The enlarged graph \(H^+\) has shifted adjacency matrix \(B^+=M_\alpha(H^+)\) and degree tuple \(\mathbf d^+\), given by
\begin{equation}
\label{eq:bord}
 B^+=\begin{pmatrix}B&-\bfb\\-\bfb^{\trans}&\alpha\end{pmatrix},
 \qquad
 \bfd^+=(d_1+b_1,\ldots,d_n+b_n,s).
\end{equation}
By \eqref{le:prin-e2}, 
the scalar Schur complement is
\begin{equation}
\label{eq:b-piv}
 q=\alpha-\bfb^{\trans}B^{-1}
 \bfb
   =\frac{\det B^+}{\det B}.
\end{equation}
A new state 
\((B^+,\bfd^+ ,H^+)\) is 
produced if $q \ge 0$. 
It is nonterminal when 
$q > 0$ and terminal
when $q = 0$. Discard it when
$q < 0$.
\end{enumerate}

By Corollary~\ref{se2-co1}
and~\eqref{le:prin-e2},
\[
q>0\iff B^+\succ0,
\qquad
q=0\iff B^+\succeq0\ \text{and }B^+\text{ is singular}.
\]

Now we define the set
\(\scrs_r(\alpha)\) of states 
and the set  \(\wisetwra\) 
of nonnegative numbers 
as follows:
\begin{enumerate}[(i)]
	\item 
\(\scrs_r(\alpha)\) is the smallest set containing the initial state \(([\alpha],(0), K_1)\) and closed under admissible updates
 from 
nonterminal states
$(B,\bfd, H)$ to states
$(B^+,\bfd^+, H^+)$
with \(q\ge 0\)
where $q$ is given in \eqref{eq:b-piv}.

\item 
 \(\wisetwra\) consists of the initial value \(\alpha\)
 and every nonnegative number
 $	q=\alpha-\bfb^{\mathsf T}B^{-1}\bfb$
 obtained from a nonterminal state $(B,\mathbf d,H)\in\mathscr S_r(\alpha)$
 and a vector $\bfb\in\{0,1\}^{|V(H)|}$ satisfying 
 \eqref{eq:admis}.
\end{enumerate}

By the definition of 
$\widehat\setw_r(\alpha)$
and $\scrs_r(\alpha)$, 
we have the following conclusion. 

\begin{lemma}\label{se2-le4}
$
0\in\widehat \setw_r(\alpha)$
if and only if $\scrs_r(\alpha)$ 
contains a terminal state.
\end{lemma}

Recall that 
$\scrq_r=\{\rho(G):G\text{ is nontrivial and connected}, \Delta(G)\le r\}$ for any
integer $r\ge 1$.
A characterization of $\scrq_r$ 
is given below.

\begin{theorem}\label{th:main}
For every \(\alpha>0\) and integer \(r\ge1\),
\begin{equation}
	\label{eq:main}
 \alpha\in \scrq_r\quad
 \Longleftrightarrow\quad
 0\in\widehat \setw_r(\alpha).
\end{equation}
There exists a nontrivial connected graph \(G\) with
\(\rho(G)=\alpha\) and \(\Delta(G)=r\) if and only if a zero output
occurs with
a
terminal update with $n + 1$ vertices must satisfy
 \(\max_{1\le i\le n+1} d_i^+=r\).
\end{theorem}

\begin{proof} ($\Leftarrow$)
Assume that $0\in\widehat \setw_r(\alpha)$. 
Then, by Lemma~\ref{se2-le4}, 
a terminal state $(B^+,\bfd^+,H^+)$
will be produced from 
a nonterminal state 
$(B, \bfd, H)\in \scrs_r(\alpha)$
by some choice 
$\bfb$ satisfying
\eqref{eq:admis}.
Thus, $q:=\alpha-\bfb^{\trans}B^{-1}
\bfb=0$,
and so Corollary~\ref{se2-co1} shows that \(M_\alpha(H^+)=B^+\) is positive semidefinite and singular. Hence \(\rho(H^+)=\alpha\) by~\eqref{eq:psd}.
Condition~\eqref{eq:admis} ensures that \(H^+\) is connected, nontrivial, and has maximum degree at most \(r\).
	
($\Rightarrow$)
Conversely, let \(G\) be nontrivial and connected with
\(\rho(G)=\alpha\) and \(\Delta(G)\le r\). 
Order its vertices \(v_1,\ldots,v_N\) so that 
$N_G(v_i)\cap \{v_1,\ldots,
v_{i-1}\}
\ne \emptyset$ for each 
$2\le i\le N$. 

For each $1\le j\le N$, 
the induced subgraph
\(G_j=G[\{v_1,\ldots,v_j\}]\) is then connected. 
For \(1\le j<N\), write
\[
1\le s=d_{G_{j+1}}(v_{j+1})
\le d_G(v_{j+1})\le r.
\]
Take
\(\bfb\) to record the neighbours of \(v_{j+1}\) in \(G_j\).
This choice satisfies \eqref{eq:admis}, because \(v_{j+1}\)
has at least one neighbour 
in $V(G_j)$ and 
\[
d_{G_j}(v_i)+b_i
=d_{G_{j+1}}(v_i)
\le d_G(v_i)\le r
\qquad(1\le i\le j).
\]
By Lemma~\ref{le:prin}, \(M_\alpha(G_j)\succ0\)
for \(j<N\). 
By \eqref{eq:psd}, 
$M_\alpha(G_N)$ is positive
semidefinite and singular.
Consequently the corresponding Schur complements 
$q$ obtained in \eqref{eq:b-piv}
are
positive until the last step, where \eqref{eq:psd} and
\eqref{eq:congr} give \(q=0\). 
By definition, $0\in\widehat \setw_r(\alpha)$.
If \(\Delta(G)=r\), the terminal
degree tuple has an entry equal to \(r\).
\end{proof}

Combining Theorem~\ref{th:main} with
 the standard bound \(\Delta(G)\le\rho(G)^2\) 
 \cite[Section 3.15, Exercise 1]
 {BrouwerHaemers}
 we obtain the following characterization.

\begin{corollary}\label{cor:all}
The spectral radii of nontrivial connected graphs satisfy
\begin{equation}\label{eq:all}
 \scrq=\{\alpha\ge1:0\in\widehat \setw_{\lfloor\alpha^2\rfloor}(\alpha)\}.
\end{equation}
\end{corollary}

The recursion has no prescribed upper bound on the number of steps.
For an algebraic target \(\alpha\), any proposed finite construction
can be checked by exact algebraic arithmetic. 
At each step in the construction, the new vertex is joined to \(s=\sum_{i=1}^{n}b_i\) vertices of the current graph \(H\). The step adds one vertex and \(s\) edges, so it increases \(|E(H)|-|V(H)|+1\) by \(s-1\). Since this quantity is initially zero for \(K_1\), 
the following observations hold:
\begin{enumerate}[(i)]
	\item the resulting graph is a tree exactly when \(s=1\) at every step, and 
	\item the resulting graph is
	 unicyclic exactly when \(s=2\) at one step and \(s=1\) at all other steps.
\end{enumerate}

\begin{figure}[htbp]
\centering
\begin{tikzpicture}[graph,>=Stealth]
  \node[vertex,label=above:{$v_1$}] (a) at (0,0) {};
  \node[vertex,label=above:{$v_1$}] (b1) at (3.15,0) {};
  \node[vertex,label=above:{$v_2$}] (b2) at (4.55,0) {};
  \node[vertex,label=left:{$v_1$}] (c1) at (7.85,-.32) {};
  \node[vertex,label=right:{$v_2$}] (c2) at (9.25,-.32) {};
  \node[vertex,label=above:{$v_3$}] (c3) at (8.55,.75) {};
  \draw (b1)--(b2);
  \draw (c1)--(c2) (c2)--(c3) (c3)--(c1);
  \draw[->] (.68,0)--(2.45,0) node[midway,above=3pt] {$q_1=\tfrac32$};
  \draw[->] (5.25,0)--(7.05,0) node[midway,above=3pt] {$q_2=0$};
  \node at (0,-.65) {$K_1$};
  \node at (3.85,-.65) {$K_2$};
  \node at (8.55,-.65) {$C_3$};
\end{tikzpicture}
\caption{The vertex additions for \(\alpha=2\) and \(r=2\) in
Example~\ref{ex:bord}.}
\label{fig:bord}
\end{figure}
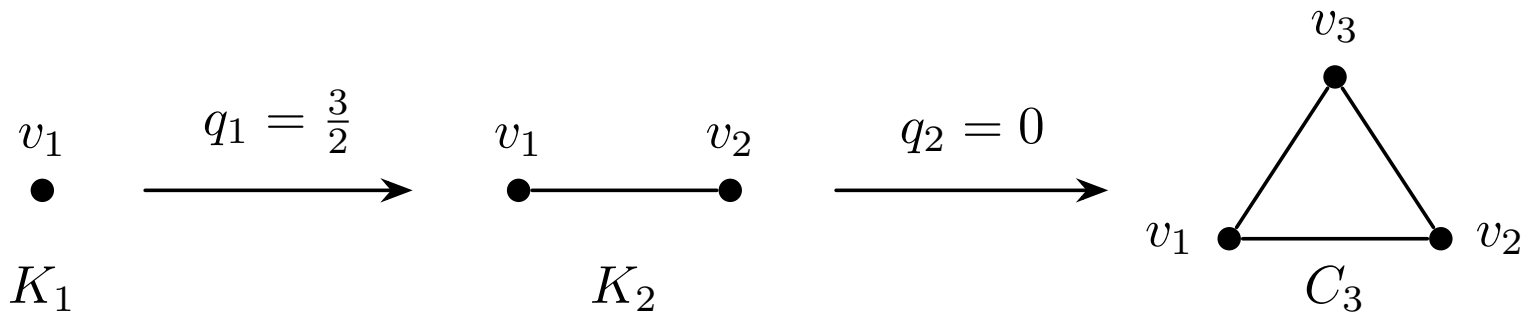

\begin{example}[A triangle]
	\label{ex:bord}
Take \(\alpha=2\) and \(r=2\). The first choice,
\(\bfb=(1)^{\trans}\), gives \(q_1=2-\tfrac12=\tfrac32\)
and produces
\[
 B=\begin{pmatrix}2&-1\\-1&2\end{pmatrix},
 \qquad B^{-1}=\frac13\begin{pmatrix}2&1\\1&2\end{pmatrix}.
\]
With \(\bfb=(1,1)^{\trans}\) in the next step,
\(\bfb^{\trans}B^{-1}\bfb=2\), so \(q_2=0\).
This produces the triangle \(C_3\), with \(\rho(C_3)=2\) and terminal
degree tuple \((2,2,2)\). Figure~\ref{fig:bord} shows both
vertex additions.
\end{example}

\section{Rooted trees and spectral reduction}	
\label{se:core}

Dong and Zhang 
\cite[Theorems 2 and 3]{DongZhang}
use a scalar recursion to characterize the existence of trees with prescribed adjacency spectral radius and maximum degree at most \(r\) or exactly \(r\). Their recursion combines the contributions of the branches at a single root. The same calculation applies when a rooted tree with a positive root label is joined to a larger graph by a single edge: eliminating the tree contributes a diagonal correction at its attachment vertex. This motivates replacing the single root by an arbitrary nonempty connected graph \(H\) and attaching rooted trees at its vertices. Theorem~\ref{se3-th1} uses the tree labels to reduce the spectral condition to a matrix indexed by \(V(H)\) and gives the corresponding maximum-degree conditions.

\subsection{The scalar recursion for trees}
\label{sec3-1}

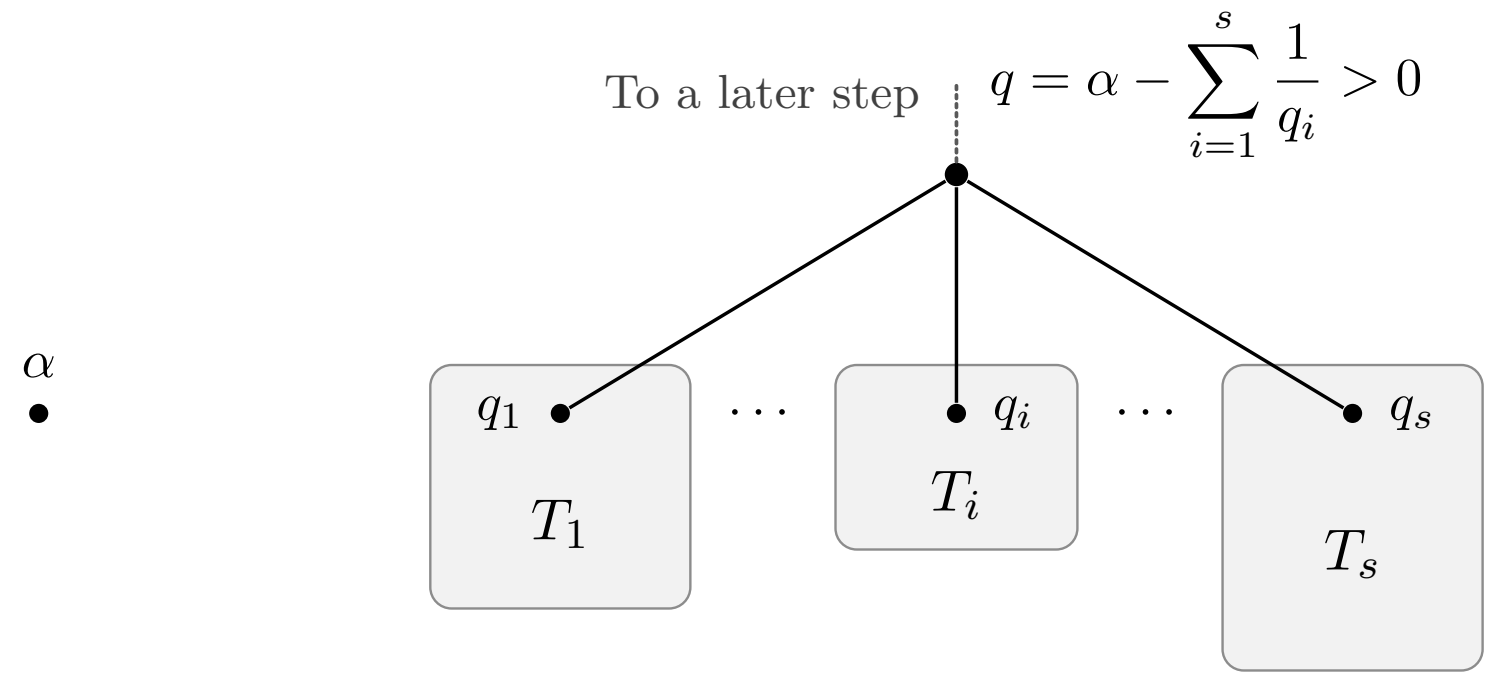
\begin{figure}[h!]
	\centering
\begin{tikzpicture}[
	font=\small,
	vertex/.style={circle, fill=black, inner sep=0pt, minimum size=3.6pt},
	tree edge/.style={line width=.65pt},
	input tree/.style={
		draw=black!45, fill=black!5, rounded corners=4pt, line width=.45pt
	},
	later edge/.style={
		draw=black!65, densely dotted, line cap=round, line width=.7pt
	}
	]
	\node[vertex] (initial) at (-1,1.28) {};
	\node[above=3pt] at (initial) {$\alpha$};
	
	\draw[input tree] (1.59,-.01) rectangle (3.31,1.60);
	\draw[input tree] (4.27,.38) rectangle (5.87,1.60);
	\draw[input tree] (6.83,-.42) rectangle (8.55,1.60);
	
	\node[vertex] (input-1) at (2.45,1.28) {};
	\node[vertex] (input-i) at (5.07,1.28) {};
	\node[vertex] (input-s) at (7.69,1.28) {};
	\node[left=3pt]  at (input-1) {$q_1$};
	\node[right=3pt] at (input-i) {$q_i$};
	\node[right=3pt] at (input-s) {$q_s$};
	
	\node[font=\normalsize] at (2.45,.55) {$T_1$};
	\node[font=\normalsize] at (5.07,.74) {$T_i$};
	\node[font=\normalsize] at (7.69,.35) {$T_s$};
	\node at (3.79,1.28) {$\cdots$};
	\node at (6.35,1.28) {$\cdots$};
	
	\node[vertex, minimum size=4.4pt] (new-root) at (5.07,2.86) {};
	\draw[tree edge] (new-root) -- (input-1)
	(new-root) -- (input-i)
	(new-root) -- (input-s);
	\node[anchor=south west, inner sep=0pt]
	at ([xshift=6pt,yshift=2pt]
	new-root.center)
	{$\displaystyle q=\alpha-\sum_{i=1}^{s}\frac{1}{q_i}>0$};
	
	\draw[later edge] (new-root) -- (5.07,3.48);
	
	\node[left=3pt, font=\footnotesize, text=black!75] at (5.07,3.38)
	{To a later step};
	
	\node at (-1,-.93) {(a) Initial tree};
	\node at (5.07,-.93) {(b) Recursion step};
\end{tikzpicture}

\caption{The rooted-tree construction for the scalar recursion. In (b), each
	shaded region represents a copy of an input tree $T_i$ with root label $q_i$
	and all leaves labelled $\alpha$. For $1\le s<r$, the $s$ copies are disjoint,
	even when input labels coincide. The trees may differ in height and
	branching pattern. The dotted stem indicates that $q$ may be used as an
	input in a later step.}
\label{f:recur-t}
\end{figure}

For \(\alpha>0\) and 
an integer \(r\ge1\), 
let \(\setw_r(\alpha)\) be the smallest
subset of \([0,\infty)\) containing \(\alpha\) and closed under the
following rule from~\cite{DongZhang}. Given positive
\(q_1,\ldots,q_s\in \setw_r(\alpha)\), with repetition allowed, form
\begin{equation}
\label{t-rule}
 q=\alpha-\sum_{i=1}^s \frac1{q_i}.
\end{equation}
Include \(q\) if \(1\le s<r\) and \(q\ge0\), or if \(s=r\) and
\(q=0\). 
Write
\(\setw_r^+(\alpha)=\setw_r(\alpha)\setminus\{0\}\).
 Every \(q\in\setw_r^+(\alpha)\) has a finite derivation from \(\alpha\) in which each application of \eqref{t-rule} uses at most \(r-1\) positive inputs and produces a positive output. Represent the initial label \(\alpha\) by a single vertex. At an application with \(s\) inputs, take one disjoint copy of a chosen rooted tree for each input,\footnote{A {\it rooted tree} is a tree in which one special vertex is singled out and designated as the root. This distinction induces a natural parent-child hierarchy and directs the edges away from the root.}
 even when input values coincide. 
 Join their roots to a new vertex and label it with the output of \eqref{t-rule}. The resulting rooted tree has root label \(q\), and every leaf has label \(\alpha\). See Figure~\ref{f:recur-t}.

Let \(\mathscr T_r(\alpha)\) be the smallest family of triples \((T,z,L)\), where \(T\) is a finite tree, \(z\) is its root, and \(L\colon V(T)\to\setw_r(\alpha)\), satisfying the following rules.
\begin{enumerate}[(i)]
	\item 
 If \(v\) is the unique vertex of \(K_1\) and \(L(v)=\alpha\), then \((K_1,v,L)\in\scrt_r(\alpha)\).

\item For \(s\in \brk{r}\), take vertex-disjoint copies \((T_i,z_i,L_i)\) of triples in \(\scrt_r(\alpha)\) with \(q_i=L_i(z_i)>0\)
for $i=1,\ldots,s$. 
Compute \(q\) by \eqref{t-rule}. 
If \(q=0\), or if \(q>0\) and \(s<r\), join a new root \(z\) to each \(z_i\) to form \(T\). Set \(L(z)=q\) and \(L|_{V(T_i)}=L_i\) for each \(i\). Then \((T,z,L)\in\scrt_r(\alpha)\).
\end{enumerate}

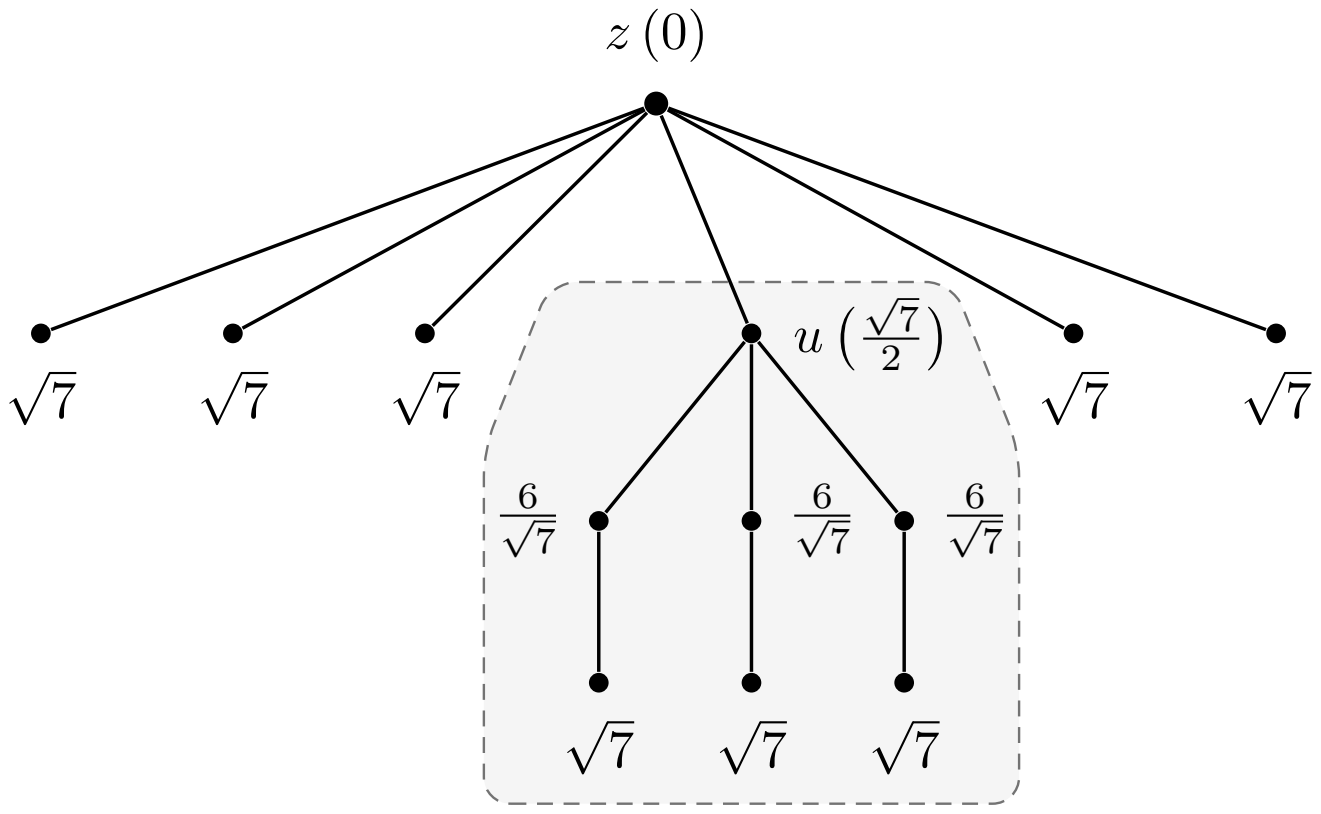
\begin{figure}[htbp]
	\centering
	\begin{tikzpicture}[
		font=\small,
		vertex/.style={circle, fill=black, inner sep=0pt, minimum size=3.7pt},
		tree edge/.style={line width=.65pt},
		subtree region/.style={
			fill=black!4, draw=black!55, densely dashed,
			rounded corners=5pt, line width=.45pt
		}
		]
		\path[subtree region]
		(-1.62,-.67) -- (1.92,-.67) -- (1.92,1.67)
		-- (1.48,2.78) -- (-1.18,2.78) -- (-1.62,1.67) -- cycle;
		
		\node[vertex, minimum size=4.5pt] (v) at (-.48,3.96) {};
		\node[vertex] (u) at (.15,2.44) {};
		\draw[tree edge] (v) -- (u);
		
		\foreach \i/\x in {1/-4.55, 2/-3.28, 3/-2.01, 4/2.28, 5/3.62} {
			\node[vertex] (leaf-\i) at (\x,2.44) {};
			\draw[tree edge] (v) -- (leaf-\i);
			\node[below=3pt] at (leaf-\i) {$\sqrt{7}$};
		}
		
		\foreach \i/\x in {1/-.86, 2/.15, 3/1.16} {
			\node[vertex] (branch-\i) at (\x,1.20) {};
			\node[vertex] (tip-\i) at (\x,.13) {};
			\draw[tree edge] (u) -- (branch-\i) -- (tip-\i);
			\node[below=3pt] at (tip-\i) {$\sqrt{7}$};
		}
		
		\node[above=4pt] at (v) {$z\,(0)$};
		\node[right=4pt] at (u) {$u\,\bigl(\frac{\sqrt{7}}{2}\bigr)$};
		\node[left=3pt]  at (branch-1) {$\frac{6}{\sqrt{7}}$};
		\node[right=3pt] at (branch-2) {$\frac{6}{\sqrt{7}}$};
		\node[right=3pt] at (branch-3) {$\frac{6}{\sqrt{7}}$};
	\end{tikzpicture}
	\caption{
	A marked member \((T,z,L)\) of \(\mathscr T_6(\sqrt7)\),
	as
	$\sqrt7-\frac5{\sqrt7}-\frac1{\sqrt7/2}=0$.
	The shaded subtree \(T_u\) is rooted at \(u\) and has labelling \(L_u=L|_{V(T_u)}\). Since \(L_u(u)=\sqrt7-\frac{3}{6/\sqrt7}=\sqrt7/2>0\), the triple \((T_u,u,L_u)\) is an unmarked member of \(\mathscr T_6(\sqrt7)\).
	}
	\label{f:rec-s7}
\end{figure}

A triple \((T,z,L)\in\mathscr T_r(\alpha)\) is called 
{\it marked} if \(L(z)=0\) and 
{\it unmarked} if \(L(z)>0\). 
Figure~\ref{f:rec-s7} shows an example of a marked triple.


Eliminate the nonroot vertices from the leaves toward the root.
When a nonroot vertex \(w\) is eliminated with a positive pivot \(q(w)\), its parent is its only remaining neighbour. By \eqref{eq:congr}, 
the elimination subtracts \(1/q(w)\) from the parent’s diagonal. Thus the pivots satisfy \eqref{t-rule}, with initial value \(\alpha\) at each leaf. For a triple \((T,z,L)\in\mathscr T_r(\alpha)\), the nonroot pivots are the positive labels \(L(w)\), and the final diagonal entry at \(z\) is \(L(z)\). Conversely, if \(M_\alpha(T)\succ0\) and every vertex has at most \(r-1\) children, all pivots are positive by \eqref{eq:congr}, and they give an unmarked triple. If such a tree is attached to a core vertex \(x\) by the edge \(xz\), eliminating \(z\) in the last step subtracts \(1/L(z)\) from the diagonal at \(x\).

By \eqref{eq:congr}, every marked triple \((T,z,L)\in\mathscr T_r(\alpha)\) has \(\rho(T)=\alpha\), since all nonroot pivots are positive and the final diagonal entry is \(L(z)=0\).
 Dong and Zhang~\cite[Theorems~2 and~3]{DongZhang}
 showed that a tree with spectral radius \(\alpha\) and maximum degree at most \(r\) exists if and only if \(0\in\mathscr W_r(\alpha)\). A tree with spectral radius \(\alpha\) and maximum degree exactly \(r\) exists if and only if \eqref{t-rule} can produce \(0\) with exactly \(r\) inputs.

\subsection{Reduction by eliminating attached trees}

Fix a nonempty connected simple graph \(H\), 
an integer 
$r\ge \max\{1,\Delta(H)\}$,
and a real number $\alpha>0$.
At each \(z\in V(H)\), choose 
an integer $s_z$ with
\(0\le s_z\le r-d_H(z)\) 
and  \(q_{z,1},\ldots,q_{z,s_z}\in \setw_r^+(\alpha)\).
Set
\begin{equation}
\label{e:br-la}
 a_z=\alpha-\sum_{\ell=1}^{s_z}\frac1{q_{z,\ell}},\qquad
 M_H(\mathbf a)=\diag(a_z:z\in V(H))-A(H).
\end{equation}
The empty sum equals zero.

Let $G$ be a connected graph, and let $H$ be a nonempty connected induced
subgraph containing every vertex of the 2-core $G'$. 
Since $V(G')\subseteq V(H)$, 
every component $T$ of
$G-V(H)$ is a tree, and 
exactly one edge $e$
joins $T$ to $H$. 
Its endpoint in $T$ is the root, and every vertex outside
$H$ has a unique parent closer to $H$.
For \(v\in V(G)\), define
\[
d_G(v,H)=\min_{u\in V(H)}d_G(v,u).
\]
Each component of \(G-V(H)\) is a tree joined to \(H\) by a single edge.
Every vertex \(v\notin V(H)\) has a unique neighbour \(p(v)\) closer to \(H\), called its parent. Its other neighbours are its children.

 For the elimination argument, assume that \(M_\alpha(T)\succ0\) for every component \(T\) of \(G-V(H)\), each rooted at its unique vertex adjacent to \(H\). At this stage no bound on \(\Delta(G)\) is imposed. Construct a labelling \(L:V(G)\to\mathbb R\) and an ordered list \(u_1,\ldots,u_t\) of the vertices outside \(H\) by the following procedure. Positive definiteness of the tree blocks ensures that every pivot at a vertex outside \(H\) is positive. The final labels at vertices of \(H\) need not be positive. When \(\Delta(G)\le r\), we will additionally verify that the resulting branch labels and attachment counts satisfy the requirements in~\eqref{e:br-la}.

\begin{enumerate}[Step (1).]
	\item 	Set \(t:=0\).
	For every \(v\in V(G)\),  set
	\[
	L(v)=\alpha,\quad
	N'(v)=
	\{u\in N_G(v)\setminus V(H):
	d_G(u,H)=d_G(v,H)+1\}.
	\]
Mark every vertex outside \(H\) as unprocessed.

\item 
If \(H=G\), terminate with \(t=0\) and an empty ordered list.

\item Choose an unprocessed vertex \(v\in V(G)\setminus V(H)\) with \(N'(v)=\varnothing\),
and set $t:=t+1$ and $u_t=v$.
 With \(w=p(v)\), update
	\[
	L(w)\leftarrow L(w)-\frac1{L(v)},
	\qquad
	N'(w)\leftarrow N'(w)\setminus\{v\},
	\]
	and mark \(v\) as processed. Retain its label \(L(v)\).

	\item Repeat Step (3) until every vertex in \(V(G)\setminus V(H)\) has been processed.
\end{enumerate}

At termination, \(t=|V(G)\setminus V(H)|\), and
the list \(u_1,\ldots,u_t\) contains every vertex outside \(H\) exactly once.
 Each child precedes its parent whenever the parent also lies outside \(H\).
Now let $G_0=G$
and $G_i=G-\{u_j: 1\le j\le i\}$
for all $i=1,\ldots,t$.
Obviously, $G_t$ is the graph $H$.

Let \(n=|V(H)|\), and order the vertices of \(G\) as \(v_1,\ldots,v_{n+t}\), where
\[
V(H)=\{v_1,\ldots,v_n\},
\qquad
v_{n+j}=u_{t-j+1}\quad
(j\in \brk{t}).
\]
Let \(A(G)\) be the adjacency matrix in this ordering, and set
\[
M_0=\alpha I_{n+t}-A(G).
\]
For \(i\in\brk{t}\), put \(r_i=n+t-i+1\). The last row and column of \(M_{i-1}\) correspond to \(v_{r_i}=u_i\). Write
\[
a_i=(M_{i-1})_{r_i r_i},
\]
and let \(q_i\) be the index such that \(v_{q_i}=p(u_i)\). Since every child precedes its parent in the elimination order, the parent is still present and \(1\le q_i<r_i\).

Let \(F_i\) be the matrix obtained 
by deleting the last row and column of \(M_{i-1}\), and define
\equ{se3-en1}
{
M_i=F_i-\frac1{a_i}e_{q_i}e_{q_i}^{\mathsf T},
}
where \(e_{q_i}\) is the coordinate unit vector in \(\mathbb R^{r_i-1}\). Thus \(M_i\) has order \(n+t-i\), and the update subtracts \(a_i^{-1}\) from the diagonal entry at the parent of \(u_i\). The preceding labelling procedure gives \(a_i=L(u_i)>0\), so each update is well-defined.

For each \(z\in V(H)\), index the components of \(G-V(H)\) joined to \(z\) as \(T_{z,1},\ldots,T_{z,s_z}\), and let \(u_{z,j}\) be the vertex of \(T_{z,j}\) adjacent to \(z\). 
Set
\[
q_{z,j}=L(u_{z,j}),\qquad a_z=L(z).
\]
Then
$
a_z=\alpha-\sum_{j=1}^{s_z}q_{z,j}^{-1}.
$
Thus the final diagonal entries have the same algebraic form as in (15). If, additionally, \(\Delta(G)\le r\), then
\[
0\le s_z=d_G(z)-d_H(z)\le r-d_H(z).
\]
Every vertex \(v\notin V(H)\) has one parent in \(G\), and therefore at most \(r-1\) children. Its positive label satisfies
\[
L(v)=\alpha-\sum_{w\text{ a child of }v}\frac1{L(w)}>0.
\]
Starting with label \(\alpha\) at vertices with no children, induction gives \(L(v)\in\mathscr W_r^+(\alpha)\). Consequently, every attached labelled rooted tree is an unmarked member of \(\mathscr T_r(\alpha)\), and all requirements in~\eqref{e:br-la} are satisfied.

By Lemma~\ref{se2-le1}, 
~\eqref{se3-en1}, and the preceding definitions, we obtain the following conclusion.


\begin{lemma}\label{se3-le1}
Let $G$ be a connected graph, and let $H$ be a nonempty
connected induced subgraph of $G$ containing every vertex
of the $2$-core $G'$. Let $\alpha>0$, and assume that
$M_\alpha(T)\succ0$ for every component $T$ of $G-V(H)$.
With the notation introduced in the preceding elimination
procedure, the following statements hold.
	\begin{enumerate}[(i)]
		\item $M_0=M_{\alpha}(G)$ and 
		$M_t=M_H(\bfa)$.
		\item 
		For every \(i\in \brk{t}\), \(M_{i-1}\) is positive definite (resp., positive semidefinite; positive semidefinite and singular) if and only if \(M_i\) has the corresponding property.
			\end{enumerate} 
\end{lemma}

\begin{proof}
Since $M_\alpha(T)\succ0$ for every attached tree $T$,
eliminating its vertices in an order in which every child
precedes its parent produces positive pivots. These pivots
are unchanged when $T$ is eliminated as part of $G$,
because the vertices of $H$ are retained throughout.
Consequently, $a_i=L(u_i)>0$ for every $i=1,\ldots,t$.

For $i\in \brk{t}$, $M_i$ is obtained as follows:
	\[
	M_{i-1}=
	\begin{pmatrix}
		F_i&-e_{q_i}\\
		-e_{q_i}^{\mathsf T}&a_i
	\end{pmatrix},
\quad 
M_i=F_i-\frac1{a_i}e_{q_i}e_{q_i}^{\mathsf T}.
	\]
After interchanging the blocks, 
	Lemma~\ref{se2-le1}
	applied to the positive scalar block \(a_i\)
	gives that 
	$
	M_{i-1}\ \text{is congruent to}\ \operatorname{diag}(a_i,M_i).
	$
	Since \(a_i>0\), the equivalences in (ii) follow. Each elimination changes only the diagonal entry at the parent. Hence, after all attached vertices have been eliminated, the remaining off-diagonal entries are those of \(-A(H)\), and the diagonal at \(z\in V(H)\) is
	\[
	L(z)=\alpha-\sum_{j=1}^{s_z}\frac1{L(u_{z,j})}.
	\]
	This proves (i), with \(q_{z,j}=L(u_{z,j})\) and \(a_z=L(z)\).
	\end{proof}

\begin{theorem} 
\label{se3-th1}
Let \(H\) be a nonempty connected simple graph, let \(\alpha>0\), and let \(r\ge\max\{1,\Delta(H)\}\) be an integer. For each \(z\in V(H)\),
let \(s_z\) be a nonnegative integer, and for \(j\in \brk{s_z}\), let \((T_{z,j},u_{z,j},L_{z,j})\in\mathcal T_r(\alpha)\) have positive root label
$
q_{z,j}=L_{z,j}(u_{z,j}).
$
Take the trees $T_{z,j}$,
where $z\in V(H)$
and \(j\in \brk{s_z}\), 
pairwise vertex-disjoint and disjoint from \(H\), and form \(G\) from $H$ and these trees
by adding the edges \(zu_{z,j}\)
for each $z\in V(H)$
and \(j\in \brk{s_z}\).
Then, the following conclusions hold:
\begin{enumerate}[(a)]
\item $\rho(G)=\alpha$ and 
$\Delta(G)\le r$ if and only if
the bounds on \(s_z\) in 
\eqref{e:br-la} hold and \(M_H(\mathbf a)\) is positive semidefinite and singular, and 

\item $\rho(G)=\alpha$ and 
$\Delta(G)=r$ if and only if 
the bounds on \(s_z\) in 
\eqref{e:br-la} hold,  \(M_H(\mathbf a)\) is positive semidefinite and singular,
and at least one of the following 
conditions holds.
\begin{enumerate}
	\item[\textup{(i)}] \(s_z=r-d_H(z)\) for some \(z\in V(H)\).
	\item[\textup{(ii)}] 
	Some attached tree \(T_{z,j}\), rooted at \(u_{z,j}\), has a vertex with exactly \(r-1\) children.
\end{enumerate}
\end{enumerate}
\end{theorem}

\begin{figure}[h!]
	\centering
	\begin{tikzpicture}
		[scale=0.85,
		font=\small,
		line cap=round,
		base graph/.style={draw=black!55,fill=black!12,line width=.5pt},
		branch/.style={draw=black!55,fill=black!7,densely dashed,line width=.5pt},
		edge/.style={draw=black,line width=.65pt},
		vertex/.style={circle,fill=black,inner sep=0pt,minimum size=4.4pt},
		]
		\path[base graph] (0,0) ellipse [x radius=2.25cm,y radius=1.30cm];
		\foreach \height in {1.75,.20,-1.75}
		\path[branch] (-4.95,\height)
		ellipse [x radius=1.15cm,y radius=.56cm];
		\foreach \height in {1.37,-1.37}
		\path[branch] (4.95,\height)
		ellipse [x radius=1.15cm,y radius=.56cm];
		
		\coordinate (z)   at (-1.67,0);
		\coordinate (w)   at (1.67,0);
		\coordinate (uz1) at (-4.00,1.75);
		\coordinate (uz2) at (-4.00,.20);
		\coordinate (uzs) at (-4.00,-1.75);
		\coordinate (uw1) at (4.00,1.37);
		\coordinate (uws) at (4.00,-1.37);
		
		\foreach \root in {uz1,uz2,uzs}
		\draw[edge] (z) -- (\root);
		\foreach \root in {uw1,uws}
		\draw[edge] (w) -- (\root);
		\foreach \point in {z,w,uz1,uz2,uzs,uw1,uws}
		\node[vertex] at (\point) {};
		
		\node at (0,.28) {$H$};
		\node[above right,inner sep=2pt] at (z) {$z$};
		\node[above left,inner sep=2pt] at (w) {$w$};
		
		\node at (-5.10,1.75) {$T_{z,1}$};
		\node at (-5.10,.20) {$T_{z,2}$};
		\node at (-4.95,-.72) {$\vdots$};
		\node at (-5.10,-1.75) {$T_{z,s_z}$};
		\node at (5.10,1.37) {$T_{w,1}$};
		\node at (4.95,0) {$\vdots$};
		\node at (5.10,-1.37) {$T_{w,s_w}$};
		
		\node[anchor=south west,inner sep=1pt] at (-3.74,1.72) {$u_{z,1}$};
		\node[anchor=south west,inner sep=1pt] at (-3.74,.29) {$u_{z,2}$};
		\node[anchor=north west,inner sep=1pt] at (-3.74,-1.72) {$u_{z,s_z}$};
		\node[anchor=south east,inner sep=1pt] at (3.74,1.34) {$u_{w,1}$};
		\node[anchor=north east,inner sep=1pt] at (3.74,-1.34) {$u_{w,s_w}$};
	\end{tikzpicture}
	\caption{
	Construction of \(G\) in Theorem~\ref{se3-th1}. The central shaded region represents \(H\), and the regions with dashed boundaries show rooted trees attached at two representative vertices \(z\) and \(w\). Each tree \(T_{x,j}\) is joined to \(x\in V(H)\) by the single edge \(xu_{x,j}\), where \(u_{x,j}\) is its root.	
}
	\label{f:attac-t}
\end{figure}
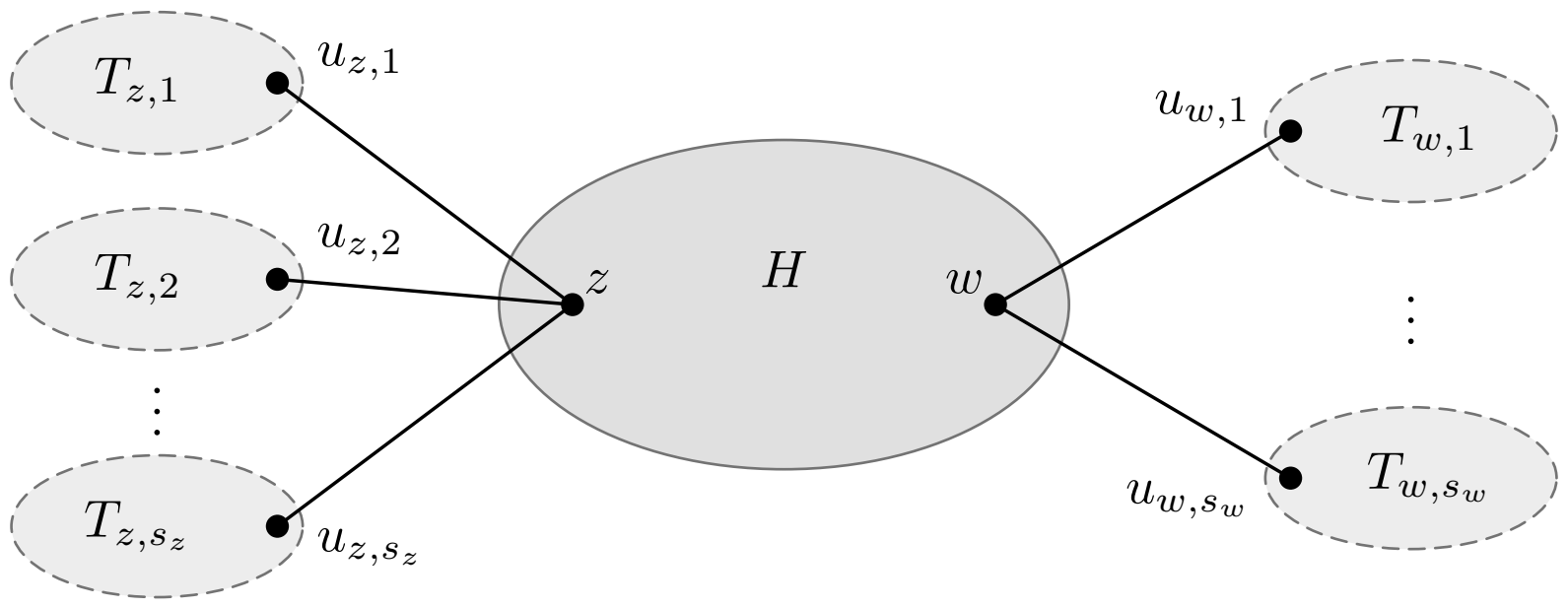

\begin{proof}
Every specified triple has positive root label. By the recursive definition of \(\mathcal T_r(\alpha)\), all labels in each attached tree are positive, and every vertex has at most \(r-1\) children.
Eliminate the vertices outside \(H\) in an order in which every child precedes its parent. By induction from the leaves and 
\eqref{t-rule}, the pivot at each \(w\in V(T_{z,j})\) equals \(L_{z,j}(w)>0\). Eliminating the root \(u_{z,j}\) therefore subtracts \(q_{z,j}^{-1}\) from the diagonal entry at \(z\). 
Lemma~\ref{se3-le1} gives
\[
M_\alpha(G)\ \text{congruent to}\
\operatorname{diag}\bigl(D,M_H(\mathbf a)\bigr),
\]
where \(D\) is the diagonal matrix of positive pivots and is omitted when there are no attached trees. Consequently, by~\eqref{eq:psd},
\[
\rho(G)=\alpha
\iff M_H(\mathbf a)\text{ is positive semidefinite and singular}.
\]
For every \(z\in V(H)\),
\[
d_G(z)=d_H(z)+s_z.
\]
Every vertex outside \(H\) has one parent and at most \(r-1\) children, so its degree is at most \(r\). Thus \(\Delta(G)\le r\) holds exactly when the bounds on \(s_z\) in \eqref{e:br-la} hold. This proves (a).

Under these bounds, a vertex of \(H\) has degree \(r\) exactly when condition (i) holds at that vertex. A vertex outside \(H\) has degree \(r\) exactly when it has \(r-1\) children, as in condition (ii). This proves (b).
\end{proof}

\noindent{\bf Remark}:
Taking \(H=K_1\) in 
Theorem~\ref{se3-th1} recovers the tree criteria in 
\cite[Theorems 2 and 3]{DongZhang}. 
If \(z\) is the sole vertex of \(H\), then \(M_H(\mathbf a)=[a_z]\), and the matrix condition reduces to
$$
a_z=\alpha-\sum_{\ell=1}^{s_z}q_{z,\ell}^{-1}=0.
$$
This is the zero-output condition in~\eqref{t-rule}, 
with \(s_z\le r\) by~\eqref{e:br-la}. For the exact-degree criterion, root the tree at a vertex \(z\) of maximum degree. Then \(\Delta(T)=r\) is equivalent to \(s_z=r\).

The same elimination gives the
strict inequality criterion.

\begin{corollary}
	\label{se3-th1-c0}
	Let \(G\), \(H\), and \(M_H(\mathbf a)\) be as in Theorem~\ref{se3-th1}, with \(|V(H)|\ge 1\). Then, 
	 \(\rho(G)<\alpha\)
	if and only if 
 \(M_H(\mathbf a)\) is positive definite.
\end{corollary} 

\begin{proof}
By~\eqref{eq:pd} and 
Lemma~\ref{se3-le1},
\[
\rho(G)<\alpha
\quad\Longleftrightarrow\quad
M_\alpha(G)\succ0
\quad\Longleftrightarrow\quad
M_H(\mathbf a)\succ0.
\]
\end{proof}



We now apply Corollary~\ref{se3-th1-c0} to proper induced subgraphs.

\begin{corollary}
\label{se3-th1-c1}
Let \(G\), \(H\), and \(M_H(\mathbf a)\) be as in Theorem~\ref{se3-th1}, with \(|V(H)|\ge2\). If \(\rho(G)\le \alpha\), 
then every nonempty proper principal submatrix of \(M_H(\mathbf a)\) is positive definite.
\end{corollary}

\begin{proof}
	Fix \(z\in V(H)\). Let \(G_z\) be obtained from \(G\) by deleting \(z\) and all vertices of the trees attached at \(z\).
	Let \(H_1,\ldots,H_k\) be the connected components of \(H-z\). For each \(i\), let \(G_i\) consist of \(H_i\) together with all trees attached to its vertices. These graphs are the connected components of \(G_z\). Since each \(G_i\) is a proper induced subgraph of the connected graph \(G\),
	Proposition 3.1.1(iii) of 
	\cite{BrouwerHaemers}
	 gives
	$
	\rho(G_i)<\rho(G)\le\alpha.
	$
	
	Corollary~\ref{se3-th1-c0} therefore gives
	$
	M_{H_i}\bigl(\mathbf a|_{V(H_i)}\bigr)\succ0
	$ for $(1\le i\le k).
	$
	
	The branch labels at the remaining vertices are unchanged. Thus, after ordering the vertices component by component, the matrix obtained from \(M_H(\mathbf a)\) by deleting the row and column corresponding to \(z\) is block diagonal with these positive definite blocks. It is therefore positive definite.
	
	Deleting further corresponding rows and columns preserves positive definiteness. Hence every nonempty proper principal submatrix of \(M_H(\mathbf a)\) is positive definite.
\end{proof}

\subsection{Graphs with a prescribed \(2\)-core}
\label{sec4}

Recall that 
\(G'\) denotes 
the \(2\)-core of \(G\), 
the graph left by
repeatedly deleting vertices of degree at most one. 

The following existence criterion follows from Theorem~\ref{se3-th1}. Indeed, the assumptions that \(H\) is connected and induced and that \(V(G')\subseteq V(H)\) ensure that \(G\) consists of \(H\) with rooted trees attached by single edges. For a realizing graph, these branches yield unmarked triples by Lemma~\ref{le:prin}
and Section~\ref{sec3-1}. Conversely, choose unmarked representatives for the prescribed root labels, using a witnessing triple when 
condition (ii) of Corollary~\ref{th:core}
applies.

\begin{corollary}
	\label{th:core}
Let \(H\) be a nonempty connected simple graph, and
let 
$r\ge \max\{1,\Delta(H)\}$
be an integer, 
  and  let 
\(\alpha>0\) be a real number.
\begin{enumerate}[(a)]
\item There exists a connected
graph \(G\) 
such that $H$ is an induced 
subgraph of $G$ with
$V(G')\subseteq V(H)$, 
 \(\Delta(G)\le r\) and
\(\rho(G)=\alpha\) if and only if there are choices in
\eqref{e:br-la} for which \(M_H(\mathbf a)\) is positive
semidefinite and singular. 
\item 
There exists a connected
graph \(G\) 
such that $H$ is an induced 
subgraph of $G$ with 
$V(G')\subseteq V(H)$, 
 \(\rho(G)=\alpha\), and \(\Delta(G)=r\) if and only if there are choices  in~\eqref{e:br-la} 
that make \(M_H(\mathbf a)\) positive semidefinite and singular and satisfy at least one of the following conditions.
\begin{enumerate}
 \item[\textup{(i)}] \(s_z=r-d_H(z)\) for some \(z\in V(H)\).
 \item[\textup{(ii)}] 
 Some \(q_{z,\ell}\) is the root label of an unmarked triple in \(\mathscr T_r(\alpha)\) whose tree has a vertex with exactly \(r-1\) children.
\end{enumerate}
\end{enumerate}
\end{corollary}

If \(H\) has minimum degree at least \(2\), then \(H\) survives every step of the 2-core pruning process. Consequently, \(V(H)\subseteq V(G')\), and the condition \(V(G')\subseteq V(H)\) gives \(G'=H\). Thus 
Corollary~\ref{th:core} yields existence criteria for graphs with a prescribed nonempty 2-core.

\section{Unicyclic graphs}
\label{se:cycle}

\subsection{Continuants for paths}

Let \(B=B(a_1,\ldots,a_n)\) denote the \(n\times n\) tridiagonal matrix 
with diagonal \(a_1,\ldots,a_n\) and adjacent off-diagonal
entries \(-1\), where \(n\ge1\).
If \(P_n=v_1\cdots v_n\) is a path, then
\[
B(a_1,\ldots,a_n)=\operatorname{diag}(a_1,\ldots,a_n)-A(P_n).
\]
The block on consecutive vertices \(v_i,\ldots,v_j\) corresponds to a subpath.

For \(1\le i\le j\le n\), let \(D(i,j)\) denote the determinant of the principal submatrix of \(B(a_1,\ldots,a_n)\) indexed by \(\llbracket i,j\rrbracket:=\{i,i+1,\ldots,j\}\). We adopt the empty-determinant convention
\[
D(i,i-1)=1\qquad(1\le i\le n+1).
\]

Expansion at the last row gives the continuant recurrence
\begin{equation}
\label{eq:cont}
 D(i,i)=a_i,\qquad
 D(i,j)=a_jD(i,j-1)-D(i,j-2)\quad(j>i).
\end{equation}

\begin{lemma}
\label{le:cont}
The matrix \(B(a_1,\ldots,a_n)\) is positive definite if and only if
\(D(1,t)>0\) for all
$t\in \brk{n}$. If \(B\) is invertible, then
\begin{equation}
\label{e:in-en}
 (B^{-1})_{11}=\frac{D(2,n)}{D(1,n)},\qquad
 (B^{-1})_{nn}=\frac{D(1,n-1)}{D(1,n)},\qquad
 (B^{-1})_{1n}=\frac1{D(1,n)}.
\end{equation}
For \(n=1\), these formulas give the same single entry \(1/a_1\).
\end{lemma}

\begin{proof}
For \(1\le t\le n\), let \(B_t=B(a_1,\ldots,a_t)\) be the leading \(t\times t\) principal submatrix of \(B\). By definition, \(\det B_t=D(1,t)\). If \(B\succ0\), then each \(B_t\succ0\), so its determinant is positive.

Conversely, suppose \(D(1,t)>0\) for every \(t\). Then \(B_1=[a_1]\succ0\). For \(t\ge2\), suppose inductively that \(B_{t-1}\succ0\). The last column of \(B_t\) above its diagonal entry is \(-e_{t-1}\), where \(e_{t-1}\) is the last coordinate unit vector in \(\mathbb R^{t-1}\). By \eqref{eq:congr}, \(B_t\) is congruent to \(\operatorname{diag}(B_{t-1},[p_t])\), where
\[
p_t=a_t-e_{t-1}^{\mathsf T}B_{t-1}^{-1}e_{t-1}.
\]
The determinant identity for the Schur complement gives
\[
p_t=\frac{\det B_t}{\det B_{t-1}}
=\frac{D(1,t)}{D(1,t-1)}>0.
\]
Thus \(B_t\succ0\), and induction gives \(B=B_n\succ0\).

Now suppose \(B\) is invertible. For \(n=1\), the empty continuants \(D(2,1)\) and \(D(1,0)\) both equal \(1\), so all three expressions in~\eqref{e:in-en} give \(1/a_1\). For \(n\ge2\), let \(B^{(j,i)}\) denote the matrix obtained from \(B\) by deleting row \(j\) and column \(i\). 
By the cofactor formula for the inverse~\cite[p.~22]
{HornJohnson}, we have
\[
(B^{-1})_{ij}
=\frac{(-1)^{i+j}\det B^{(j,i)}}{\det B}.
\]
Deleting row and column \(1\) leaves the principal block on indices \(2,\ldots,n\), whose determinant is \(D(2,n)\). Deleting row and column \(n\) leaves \(B_{n-1}\), whose determinant is \(D(1,n-1)\). This gives the two diagonal entries in~\eqref{e:in-en}.

For \((B^{-1})_{1n}\), the cofactor uses row \(n\) and column \(1\) because the adjugate is the transpose of the cofactor matrix. Set \(C=B^{(n,1)}\). 
Note that 
\(C\) is lower triangular with determinant \((-1)^{n-1}\). Thus,
\[
(B^{-1})_{1n}
=\frac{(-1)^{n+1}(-1)^{n-1}}{D(1,n)}
=\frac{1}{D(1,n)}.
\]
\end{proof}

\subsection{Spectral-radius criteria}

Fix integers \(m\ge3\) and \(r\ge2\), and a real number \(\alpha>0\). Label the vertices of \(C_m\) in cyclic order as \(v_1,\ldots,v_m\), and choose the branch labels in \eqref{e:br-la} with \(H=C_m\).
 Write \(s_i=s_{v_i}\)
and \(a_i=a_{v_i}\), so \(0\le s_i\le r-2\).
Use \(D(i,j)\) for the continuants of the path on
\(v_1,\ldots,v_{m-1}\) with diagonal \(a_1,\ldots,a_{m-1}\).
Thus, by \eqref{eq:cont}, 
\begin{equation}
	\label{eq:contin}
	\begin{aligned}
&D(i,i-1)=1\, (1\leq i\leq m),
\qquad
&&D(i,i)=a_i\, (1\leq i\leq m-1),\\
&D(i,j)=a_jD(i,j-1)-D(i,j-2),
		\qquad &&(1\leq i<j\leq m-1).
	\end{aligned}
\end{equation}

\begin{theorem}
\label{th:cycle}
There exists a unicyclic graph \(G\) of cycle length \(m\) with
\(\Delta(G)\le r\) and \(\rho(G)=\alpha\) if and only if the
choices in \eqref{e:br-la} can be made so that
\begin{equation}
\label{e:cyc-t}
 D(1,t)>0\quad(1\le t\le m-1),\qquad
 a_m D(1,m-1)=D(1,m-2)+D(2,m-1)+2.
\end{equation}
There exists such a graph with \(\Delta(G)=r\) if and only if
the choices also have \(s_i=r-2\) for some \(i\), or
some chosen branch label
\(q_{v_i,\ell}\), 
where $1\le i\le m$
and $1\le \ell\le s_i$, 
is the root label of an unmarked triple in \(\mathscr T_r(\alpha)\) whose tree has a vertex with exactly \(r-1\) children.
\end{theorem}

\begin{proof}
We may apply Corollary~\ref{th:core} with \(H=C_m\).
For the choices in \eqref{e:br-la}, write \(a_i=a_{v_i}\) and \(M=M_{C_m}(\mathbf a)\). The vertices \(v_1,\ldots,v_{m-1}\) induce a path, so the principal submatrix of \(M\) on these vertices is \(B=B(a_1,\ldots,a_{m-1})\). The two neighbours of \(v_m\) on the cycle are \(v_1\) and \(v_{m-1}\). Hence, with \(v_m\) placed last,
\equ{th4.2-e1}
{
M=\begin{pmatrix}
	B&-h\\
	-h^{\mathsf T}&a_m
\end{pmatrix},
\qquad h=e_1+e_{m-1},
}
where the \(e_i\) are the coordinate unit vectors in \(\mathbb R^{m-1}\).
Eliminating the attached trees gives the positive semidefinite singular matrix \(M\).

By Corollary~\ref{se3-th1-c1},
\(B\succ0\). 
Lemma 4.1 now gives \(D(1,t)>0\) for \(1\le t\le m-1\).
By~\eqref{eq:congr}, \(M\) is congruent to 
$ \diag\bigl(B,
[a_m-h^{\mathsf T}B^{-1}h]
\bigr)$.
Using the symmetry of \(B^{-1}\) and applying~\eqref{e:in-en} with \(n=m-1\), we obtain
\equ{e:cyc-sc}
{
\begin{aligned}
	a_m-h^{\mathsf T}B^{-1}h
	&=a_m-\bigl((B^{-1})_{11}
	+(B^{-1})_{m-1,m-1}
	+2(B^{-1})_{1,m-1}\bigr)\\
	&=a_m-\frac{D(2,m-1)+D(1,m-2)+2}{D(1,m-1)}.
\end{aligned}
}
Since \(M\) is singular and \(B\) is invertible, the scalar in \eqref{e:cyc-sc} is zero. Multiplying by \(D(1,m-1)>0\) gives the equality in~\eqref{e:cyc-t}.

Conversely, suppose the choices in \eqref{e:br-la} satisfy \eqref{e:cyc-t}. The inequalities in \eqref{e:cyc-t} give \(B\succ0\) by Lemma~\ref{le:cont}, and the equality makes the scalar in \eqref{e:cyc-sc} zero. Thus \(M\) is congruent to \(\operatorname{diag}(B,[0])\), so it is positive semidefinite and singular. Corollary~\ref{th:core}(a) gives a connected graph \(G\) with \(2\)-core \(C_m\), \(\Delta(G)\le r\), and \(\rho(G)=\alpha\). It follows that \(G\) is unicyclic with cycle length \(m\).

For the exact degree assertion, a cycle vertex \(v_i\) has degree \(2+s_i\), which equals \(r\) precisely when \(s_i=r-2\). A vertex outside the cycle has one parent in \(G\), so it has degree \(r\) precisely when it has \(r-1\) children in its attached rooted tree. These are the two conditions in Corollary~\ref{th:core}(b) specialized to \(H=C_m\).	
\end{proof}

\begin{example}[A degree attained in an attached tree]
\label{ex:cycle-mark}
Let \(G\) be obtained by joining \(v_3\) of \(C_3\) to the centre \(z\) of a disjoint \(K_{1,3}\), as in Figure 5. Then \(\Delta(G)=4\), attained at \(z\), whereas the largest degree on the cycle is \(3\).

Apply Theorem~\ref{th:cycle} with \(m=3\) and \(r=4\). For \(\alpha>2\), the three leaves have label \(\alpha\), so the branch root has positive label
$
q=\alpha-\frac3\alpha.
$
Thus
\[
s_1=s_2=0,\qquad s_3=1,\qquad
a_1=a_2=\alpha,\qquad a_3=\alpha-\frac1q.
\]
The required continuants are
\[
D(1,1)=D(2,2)=\alpha,\qquad
D(1,2)=\alpha^2-1.
\]
	Thus, the conditions in 
	\eqref{e:cyc-t} give 
	the equation below
	\begin{align}
	a_3D(1,2)
	&=D(1,1)+D(2,2)+2;
	\nonumber
	\\
	\left(\alpha-\frac1{\alpha-3/\alpha}\right)(\alpha^2-1)
	&=2\alpha+2;
		\nonumber \\
		\label{e:tri-s}
		\alpha-\frac{1}{\alpha-3/\alpha}&=\frac{2}{\alpha-1}.
	\end{align}
The largest real root of the above equation is 
$\alpha\approx 2.368649372$.
	Hence the spectral radius is  approximately
	\(2.36865\).
\end{example}

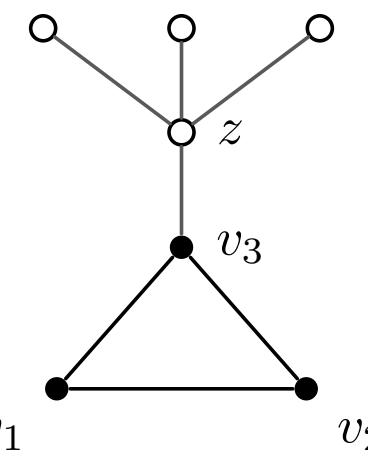
\begin{figure}[h!]
\centering
\begin{tikzpicture}
[graph,y=.8cm]

  \begin{scope}[xshift=5cm]
    \node[vertex,label=below left:{$v_1$}] (b1) at (0,0) {};
    \node[vertex,label=below right:{$v_2$}] (b2) at (1.65,0) {};
    \node[vertex,label=right:{$v_3$}] (b3) at (.825,1.17) {};
    \node[branchvertex,label=right:{$z$}] (b4) at (.825,2.12) {};
    \node[branchvertex] (b5) at (-.09,2.98) {};
    \node[branchvertex] (b6) at (.825,2.98) {};
    \node[branchvertex] (b7) at (1.74,2.98) {};
    \draw (b1)--(b2) (b2)--(b3) (b3)--(b1);
    \draw[branch] (b3)--(b4) (b4)--(b5) (b4)--(b6) (b4)--(b7);
  \end{scope}
\end{tikzpicture}
\caption{
A unicyclic graph with core \(C_3\) and maximum degree \(4\), attained at the branch root \(z\). Open vertices belong to the attached tree.
}
\label{fi:uni}
\end{figure}

\section{Open problems}
\label{sec7}

 \begin{prob}
[Effective recognition with a prescribed \(2\)-core] 
\label{prob1}
Let \(H\) be a connected graph of minimum degree at least two, let \(r\ge\Delta(H)\) be an integer, and let \(\alpha>0\) be a real algebraic number specified by its minimal polynomial and a rational isolating interval. Is there an algorithm that always decides whether there exists a connected graph \(G\) with
\[
G'\cong H,\qquad \rho(G)=\alpha,\qquad \Delta(G)\le r?
\]
Does such an algorithm exist for \(\Delta(G)=r\)? In each case, an equivalent question is whether one can compute an order bound for every input \((H,r,\alpha)\) that bounds the order of a smallest realizing graph whenever one exists. A concrete first case is \(H=C_4\) and \(r=3\).
\end{prob}


 \begin{prob}
[Square-root spectral radii of unicyclic graphs]
	\label{prob2}
For a nonsquare integer \(k\ge5\) and an even integer \(m\ge4\), determine all integers \(r\) for which a unicyclic graph of cycle length \(m\) has spectral radius \(\sqrt{k}\) and maximum degree \(r\). 
 \end{prob}
 
 Dong and Zhang~\cite{DongZhang} settle the corresponding problem for trees in \cite{DongZhang}. 

If $k$ is nonsquare, the minimal polynomial of $\sqrt{k}$ over $\mathbb Q$
is $x^2-k$. Since this polynomial divides $\det(xI-A(G))$, the number
$-\sqrt{k}$ is also an adjacency eigenvalue. For a connected graph with
$\rho(G)=\sqrt{k}$,
Proposition 3.4.1(ii) of \cite{BrouwerHaemers} implies
that $G$ is bipartite. Its unique cycle therefore has even length.
Moreover, any realizing graph satisfies
\begin{equation}
\label{eq:problem2}
	\lceil\sqrt{k}\rceil\leq r\leq k-1.
\end{equation}
Indeed, $\rho(G)\leq\Delta(G)=r$ gives the lower bound. The star formed
by the edges incident with a degree-$r$ vertex is a proper subgraph of
the connected unicyclic graph $G$, so strict spectral monotonicity gives
$\sqrt r<\rho(G)=\sqrt k$, proving the upper bound.

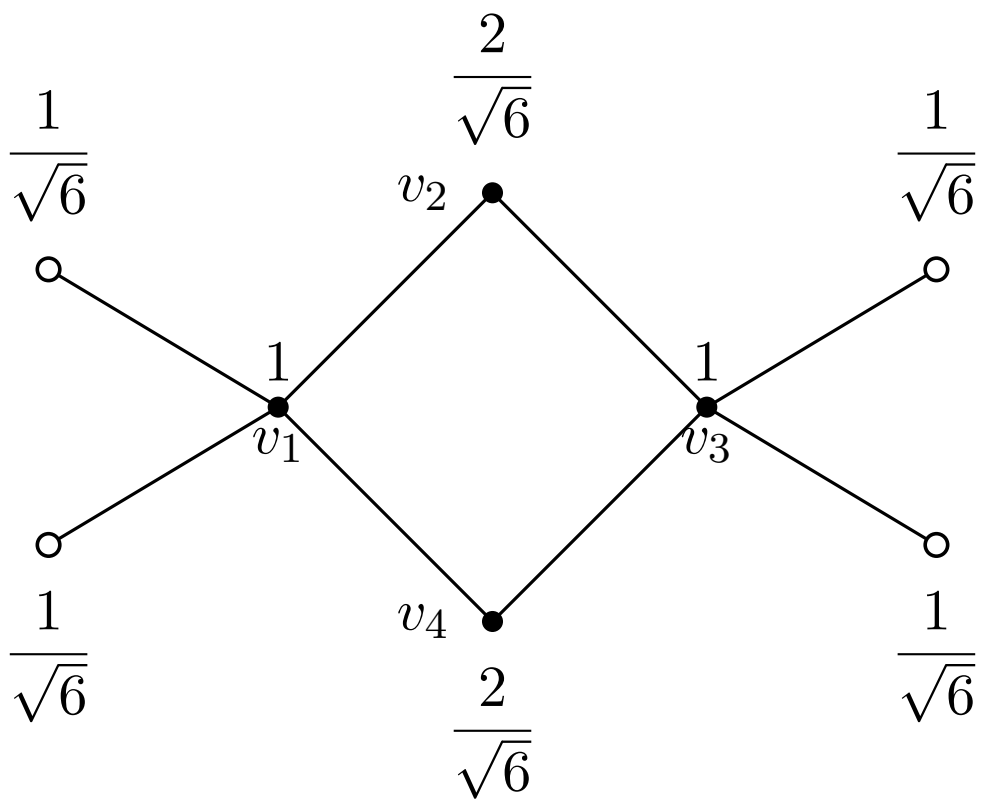
\begin{figure}[h!]
\begin{center}
	\begin{tikzpicture}[scale=0.9,
		x=1cm,y=1cm,
		line width=0.60pt,
		line cap=round,
		line join=round,
		every node/.style={font=\normalsize,inner sep=0pt,outer sep=0pt}
		]
		\def\FilledRadius{2.20pt}
		\def\OpenRadius{2.35pt}
		
		\coordinate (v1) at (-1.5750, 0);
		\coordinate (v2) at ( 0,      1.5750);
		\coordinate (v3) at ( 1.5750, 0);
		\coordinate (v4) at ( 0,     -1.5750);
		
		\coordinate (lu) at (-3.2625, 1.0125);
		\coordinate (ld) at (-3.2625,-1.0125);
		\coordinate (ru) at ( 3.2625, 1.0125);
		\coordinate (rd) at ( 3.2625,-1.0125);
		
		\draw (v1)--(v2)--(v3)--(v4)--cycle;
		\draw (lu)--(v1)--(ld);
		\draw (ru)--(v3)--(rd);
		
		\foreach \p in {v1,v2,v3,v4}
		\fill[black] (\p) circle[radius=\FilledRadius];
		
		\foreach \p in {lu,ld,ru,rd}
		\draw[draw=black,fill=white,line width=0.68pt]
		(\p) circle[radius=\OpenRadius];
		
		\node[below=4pt] at (v1) {$v_1$};
		\node[left=8pt]  at (v2) {$v_2$};
		\node[below=4pt] at (v3) {$v_3$};
		\node[left=8pt]  at (v4) {$v_4$};
		
		\node[above=5pt] at (v1) {$1$};
		\node[above=5pt] at (v3) {$1$};
		\node[above=9pt] at (v2) {$\dfrac{2}{\sqrt{6}}$};
		\node[below=9pt] at (v4) {$\dfrac{2}{\sqrt{6}}$};
		
		\node[above=9pt] at (lu) {$\dfrac{1}{\sqrt{6}}$};
		\node[below=9pt] at (ld) {$\dfrac{1}{\sqrt{6}}$};
		\node[above=9pt] at (ru) {$\dfrac{1}{\sqrt{6}}$};
		\node[below=9pt] at (rd) {$\dfrac{1}{\sqrt{6}}$};
	\end{tikzpicture}%
\end{center}

\caption{A unicyclic graph with 
core $C_4$, maximum degree $4$ and spectral radius 
$\sqrt 6$.}
\label{f7}
\end{figure}

The cycle length must be even because \(-\sqrt{k}\) is an algebraic conjugate of \(\sqrt{k}\).
Since $G$ is connected and $-\rho(G)$ is an adjacency
eigenvalue, $G$ is bipartite by
\cite[Proposition~3.4.1(ii)]{BrouwerHaemers}.
One realizable value is \(r=k-2\) for every even \(m\). Attach \(k-4\) leaves to each alternate vertex of \(C_m\). Assign values \(1\), \(2/\sqrt{k}\), and \(1/\sqrt{k}\) to the decorated cycle vertices, the other cycle vertices, and the leaves, respectively. These values form a positive eigenvector for \(\sqrt{k}\). 
An example is shown in Figure~\ref{f7} with $m=4$ and
$k=6$.

\begin{prob}
	[Infinitely many realizations of one spectral radius]
	\label{prob3}
	For fixed \(H\) and \(r\) as in Problem~\ref{prob1}, characterize the values of \(\alpha\) for which there are infinitely many pairwise nonisomorphic connected graphs \(G\) with
\[
G'\cong H,\qquad \rho(G)=\alpha,\qquad \Delta(G)=r.
\]
\end{prob}
Seeger and Sossa~\cite{Seeger}
construct infinite families with a common spectral radius by varying the length of a cycle with periodically attached graphs. Here the entire \(2\)-core is prescribed.

\begin{prob}
[Accumulation points for a prescribed cycle]
\label{prob4} 
For fixed integers \(m\ge3\) and \(r\ge3\), determine the accumulation points of
\[
\{\rho(G): G\text{ is connected},\ G'\cong C_m,\ \Delta(G)\le r\}.
\]
\end{prob}
This is a fixed-cycle version of Hoffman’s limit-point problem
\cite{Wang}.
Let $C_m^+$ be $C_m$ with one pendant vertex. Every graph with 2-core
$C_m$ is either $C_m$ itself or contains $C_m^+$ as a subgraph. Consequently
\begin{equation}\label{eq:fixed-cycle-gap}
	\rho(G)=2\quad\text{if }G=C_m,
	\qquad
	\rho(G)\geq\rho(C_m^+)>2\quad\text{otherwise}.
\end{equation}
Thus $2$ is isolated in the set in Problem 4, unlike its role as a limit
point when cycle lengths are allowed to vary.

\sgap
\noindent\textbf{Declaration on the use of AI:}
During the preparation of this work, the author used AI systems to assist in exploring a
proof of Theorem~\ref{se3-th1}. All AI-generated suggestions were verified and refined by the author,
who takes full responsibility for the correctness and originality of the paper.

\sgap

\noindent{\bf Declaration of competing interest}

There is no conflict of interest.\\

\noindent{\bf Data availability}

No datasets were
generated or analysed
during this study.




\end{document}